\documentclass[a4paper,11pt,english]{scrartcl}
\usepackage{pdfsync}

\KOMAoptions{headinclude=true,footinclude=false}
\KOMAoptions{DIV=12}
\KOMAoptions{parskip=half+}
\usepackage[includehead=true,hmargin=3cm,vmargin=3cm]{geometry}

\usepackage[utf8]{inputenc}      
\usepackage[automark]{scrpage2}  
\usepackage{abdo}                

 \setheadsepline{1pt}
  \ohead{\pagemark} \ihead{\headmark}

\newcommand{\Mpt}{\caM^*(\mrm{pt})}
\newcommand{\LHpt}{LH(\mrm{pt})}
\newcommand{\suhi}{\mrm{hi}}

\begin{document}

\title{Morphic cohomology of toric varieties}
\date{January 18, 2010}
\begin{authors}
  \authoritem{Abdó Roig-Maranges}
             {abdo.roig@upc.edu} 
             {} 
             {Partially supported by MICINN MTM2009-09557 from the
               Spanish government and the ``Comissionat per a
               Universitats i Recerca'' from the ``Generalitat de
               Catalunya''.}

   \affiliation{Universitat Politècnica de Catalunya\\
               Dept. de Matemàtica Aplicada 1}
              {Av. Diagonal, 647\\
               08028 Barcelona (Spain)}
\end{authors}
\maketitle

\begin{abstract}
  In this paper we construct a spectral sequence computing a modified
  version of morphic cohomology of a toric variety (even when it is
  singular) in terms of combinatorial data coming from the fan of the
  toric variety.
\end{abstract}

\tableofcontents

\section{Introduction}
\label{sec:intro}

Morphic cohomology is a cohomological theory on algebraic varieties
introduced by Friedlander and Lawson in \cite{FriLaw:Cocycles}. On one
hand this theory has a very geometric definition, and on the other it
is strongly related to the abstractly defined motivic cohomology
theory.

In this paper, we define, along the lines of \cite{Fri:BlochOgus}, a
modification of morphic cohomology in order to have Mayer-Vietoris and
homotopy invariance properties of the theory for all quasi-projective
varieties (not only smooth). Then, we use those properties to
construct a spectral sequence computing the morphic cohomology of a
toric variety in terms of its combinatorial data.

In section \ref{sec:morphic} we recall the definition of morphic
cohomology and the theorems we need, and define the modification of
the theory we will use. This is mainly an expository section except,
maybe, subsection \ref{sec:mormod}.

In section \ref{sec:toric} we recall the definition of toric variety,
and the constructions we will use, mainly to set up the notation. We
also write down an explicit computation of the morphic cohomology of
algebraic tori.

Finally, in section \ref{sec:sseq} we build a resolution of the
constant sheaf $\bbZ_X$ on a toric variety in terms of the
combinatorial data. This resolution allows us to construct a spectral
sequence (theorem \ref{thm:sseq}) computing the hypercohomology of a
complex of sheaves $\caF^*$ on the toric variety $X(\Delta)$ in terms
of the combinatorial data and the value of this hypercohomology on
algebraic tori. Then we specialize this to the case of morphic
cohomology, giving a very explicit spectral sequence converging to the
morphic cohomology of $X(\Delta)$ and whose second page involves only
combinatorial data. Moreover, we prove its rational degeneration
(theorem \ref{thm:sseqrat}). Finally, we also use the spectral
sequence to extend the Suslin conjecture, proved in
\cite{FriHaWa:TechKth} and \cite{Voi:SemiKth} for smooth linear
varieties, to all toric varieties, even singular ones.

I want to express my gratitude to my advisor, Pere Pascual, for his
help and encouragement during this work, and also to Eric Friedlander
and Mircea Voineagu, from whom I learnt a lot during the time I spent
in Los Angeles.

\section{Morphic cohomology}
\label{sec:morphic}

In this section we recall some facts about morphic cohomology. Let $Y$
be a projective variety over $\bbC$, with a fixed projective
embedding. It is a classical fact (see chapter 1 in
\cite{Kollar:RatCur}) that the set of effective $k$-cycles on $Y$ has
the structure of an algebraic variety. This variety is called the
\emph{Chow variety}, and we will denote it by $C(Y,k)$. It may have
infinitely many connected components, corresponding to the homology
classes of the cycles. The Chow variety has a natural operation given
by the sum of cycles. This operation is algebraic, so $C(Y,k)$ is a
monoid in the category $\ctProj_{\bbC}$ of projective varieties (see
\cite{Fri:Lawson}).

Set theoretically, the group of algebraic $k$-cycles on $Y$ is the
group completion of the Chow variety,
\begin{equation}
  Z(Y,k) = C(Y,k)^+ = C(Y,k)\tms C(Y,k) /
  C(Y,k).
  \label{eq:defcyc}
\end{equation}
This can be used to induce a topology on $Z(Y,k)$ from the
analytic topology on the Chow variety. As shown in \cite{Law:Cycles}
and \cite{Fri:Lawson}, this topology on $Z(Y,k)$ does not depend
on the projective embbeding of $Y$. Moreover, its homotopy type is a
very interesting invariant containing information about the geometry
of $Y$. The homotopy groups of this space $Z(Y,k)$ are called
the \emph{Lawson homology groups} of $Y$, and are denoted by
$L_kH_n(Y)$. The usual indexing convention is as follows,
\begin{equation}
  L_kH_n(Y) = \pi_{n-2k}Z(Y,k).
  \label{eq:deflaw}
\end{equation}

The construction of the Chow varieties $C(Y,k)$ works only for
projective $Y$, but the topology on $Z(Y,k)$ can be defined for
quasi-projective varieties (see \cite{Lim:Qproj}). Take $Y\subset
\overline{Y}$ a projective closure of the quasi-projective variety
$Y$, and let $Y^{\oo} = \overline{Y}\setm Y$, which is again
projective. Then we define the topology in $Z(Y,k)$ as the quotient
topology given by the set-theoretic identification
\begin{equation}
  Z(Y,k) = Z(\overline{Y},k)/Z(Y^{\oo},k).
  \label{eq:defcycqproj}
\end{equation}
Theorem 4.3 in \cite{Lim:Qproj} states that the topology on $Z(Y,k)$
we have just described does not depend on the choice of projective
closure $\overline{Y}$.

It will be useful to work simplicially, so we denote the singular
simplicial set associated to $Z(Y,k)$ by
\begin{equation*}
  Z_{\bt}(Y,k) = \mrm{Sing}_{\bt}Z(Y,k).
\end{equation*}
The group law on $Z(Y,k)$ makes $Z_{\bt}(Y,k)$ into a simplicial
abelian group.

One important fact about those spaces of cycles is that the homotopy
groups of the spaces of 0-cycles coincide with the Borel-Moore
homology groups of $Y$ (singular homology if $Y$ is projective),
\begin{equation}
  \pi_lZ(Y,0) \simeq H^{\mrm{BM}}_l(Y,\bbZ).
  \label{eq:DT}
\end{equation}
This is a consequence of a classical result of Dold and Thom
\cite{DoldThom:Symmetric}. From this we see that the space of cycles
$Z(\bbA^q,0)$ has only nontrivial homotopy groups in dimension
$2q$, so it is an algebraic model of an Eilenberg-Maclane space
$K(\bbZ, 2q)$.

This interpretation motivates the construction of a cohomological
theory related to the Lawson homology groups, \emph{morphic
  cohomology}, in the same way singular cohomology is related to
singular homology. The inspiration for constructing morphic cohomology
comes from the fact that singular cohomology can be interpreted as
homotopy groups of spaces of continuous maps to an Eilenberg-Maclane
space,
\begin{equation}
  H^{2q-l}(X, \bbZ) \simeq \pi_l \Map(X, K(\bbZ, 2q)).
  \label{eq:cohmaps}
\end{equation}
Then one could think of the homotopy groups of spaces of ``algebraic
maps'' from an algebraic variety $X$ to our algebraic model of
$K(\bbZ, 2q)$, the space of 0-cycles $Z(\bbA^q, 0)$. This is roughly
the idea behind the construction of morphic cohomology, as we will
recall in section \ref{sec:morcpx}. But before this, it will be
convenient to talk about a more general bivariant theory depending on
two algebraic varieties, $X$ and $Y$.

\subsection{Bivariant Morphic spaces}

We are going to describe a bivariant theory, following
\cite{Fri:CocyclesQproj}. This theory depends on two algebraic
varieties $X$ and $Y$, and roughly speaking, is defined as the
homotopy groups of spaces of ``algebraic maps'' from $X$ to a space of
cycles $Z(Y,k)$. The first thing to do is to make precise what we mean
by ``algebraic maps''. For technical reasons, it happens to be more
convenient to allow a slightly more general class of maps than the
usual algebraic morphisms of varieties.

\begin{dfn} 
  \label{dfn:contalg}
  A \emph{continuous algebraic map} between two algebraic varieties
  $X$ and $Y$ is a rational map $f\cl X \rto Y$ which is defined at
  every point and is continuous. We will denote by $\Mor(X,Y)$ the set
  of continuous algebraic maps.
\end{dfn}

\begin{rmk}
  If $X$ is a normal variety, a continuous algebraic map $f\cl X \to
  Y$ is a morphism. This is a consequence of Zariski's main theorem
  applied to the projection $\Gamma_f \to X$ of the graph of the map
  to $X$.

  In particular, this means that the technical point involving
  continuous algebraic maps appears only for non-normal varieties.

  For more details on the role of continuous algebraic maps see
  \cite{Fri:BlochOgus} section 1 and \cite{FriWa:Function}.
\end{rmk}

We want to topologize the spaces $\Mor(X,Y)$. The first choice of
topology would be the compact-open topology, however it is not the
right one. Its problems come from two different sources: one appears
when $X$ is non-complete, and the other when $X$ is non-normal. Here
we just describe the correct topology we need, and for the
justification of why this seems to be the right choice we refer to
\cite{FriLaw:Duality} appendix C (for the case in which $X$ is normal)
and \cite{Fri:BlochOgus} proposition 1.4 for the general case.

\begin{dfn}
  \label{dfn:mortop}
  Let $X$ be a quasi-projective variety and $Y$ a projective one.

  When $X$ is normal, we provide the space $\Mor(X,Y)$ of continuous
  algebraic maps with a topology such that a sequence of maps
  $\set{f_n}$ converges to $f \in \Mor(X,Y)$ if, and only if
  \begin{enumerate}
  \item it converges for the compact-open topology,
  \item there is a compactification $X \subset \overline{X}$ for which
    the sequence of closures of graphs $\overline{\Gamma_{f_n}}
    \subset \overline{X}\times Y$ has bounded degree.
  \end{enumerate}

  When $X$ is not normal, we provide the space $\Mor(X,Y)$ with the
  topology induced as a subspace of $\Mor(\tilde{X}, Y)$, where
  $\tilde{X}\to X$ is the normalization of $X$.
\end{dfn}

\begin{rmk}
  In the case in which $X$ is projective and normal, the condition on
  the boundedness of the degree becomes void, and then the topology on
  $\Mor(X,Y)$ coincides with the compact-open topology.
\end{rmk}

There is one last technical point we must deal with. In the situation
in which we are interested, $\Mor(X, Z(Y,k))$, the second argument is
the space of cycles $Z(Y,k)$ which is not an algebraic variety. But it
is the group completion of the Chow variety $C(Y,k)$, so we can
consider instead the well defined space $\Mor(X, C(Y,k))^+$ where
$(\hy)^+$ denotes a group completion as in \eqref{eq:defcyc}.

Moreover, if $Y$ is quasi-projective, we can perform a similar trick
as for $Z(Y,k)$, but on the whole $\Mor(X, C(Y,k))^+$. We chose a
projective embedding for $Y$, and take the projective closure
$\overline{Y}$, with closed complement $Y^{\oo}$. In this case though,
the bivariant theory will depend on the pair $(\overline{Y},
Y^{\oo})$, not only on $Y$. Finally, as it will be convenient to work
simplicially, we take singular chains. This is the content of the
following definition.
\begin{dfn}
  \label{dfn:morsp}
  Let $X$ be quasi-projective variety and $Y$ a projective one. Then,
  the \emph{bivariant morphic space} $M_{\bt}(X,Y,k)$ is the
  simplicial abelian group
  \begin{equation}
    M_{\bt}(X,Y,k) = (\Sing_{\bt}\Mor(X, C(Y,k)))^+,
  \end{equation}
  where the group completion $(\hy)^+$ here means a levelwise
  group-completion of the simplicial abelian monoid.

  If $Y$ is quasi-projective, with projective closure $\overline{Y}$
  and closed complement $Y^{\oo}$, Then the \emph{bivariant morphic
    space} is the quotient simplicial abelian group
  \begin{equation}
    M_{\bt}(X,\overline{Y}/Y^{\oo},k) = M_{\bt}(X,\overline{Y},k) / M_{\bt}(X,Y^{\oo},k).
  \end{equation}
\end{dfn}

\begin{rmk}
  The construction of this morphic space is contravariantly functorial
  in $X$, by composition. It is also covariantly functorial for proper
  morphisms of the pair $(\overline{Y},Y^{\oo})$. That is, if
  $\overline{Y}\to \overline{Y}'$ is a proper morphism sending
  $Y^{\oo}$ to $Y'^{\oo}$, the induced map $C(\overline{Y},k) \to
  C(\overline{Y}',k)$ sends the cycles in $C(Y^{\oo},k)$ to
  $C(Y'^{\oo},k)$ so we get a continuous map
  \begin{equation*}
    \xymatrix{M_{\bt}(X,\overline{Y}/Y^{\oo},k) \ar[r] & M_{\bt}(X,\overline{Y}'/Y'^{\oo},k).}
  \end{equation*}
\end{rmk}

Those spaces have a localization propery with respect to the second
argument analogous to proposition 3.2 in \cite{Lim:Qproj}.
\begin{thm}
  \label{thm:loc}
  Let $X$ be a quasi-projective variety and let $Y_0 \subset Y_1
  \subset Y_2$ a triple of projective varieties. Then, the sequence of
  simplicial abelian groups
  \begin{equation}
    M_{\bt}(X, Y_1/Y_0,k) \to M_{\bt}(X, Y_2/Y_0, k) \to M_{\bt}(X,
    Y_2/Y_1, k)
  \end{equation}
  is a fibration sequence and, in particular, there is an associated
  long exact sequence of homotopy groups. 
\end{thm}
\begin{prf}
  The proof is analogous as proposition 3.2 in \cite{Lim:Qproj} using
  the fact that a quotient of simplicial abelian groups produces a
  fibration sequence (see \cite{GoeJar:Simplicial} corollary V.2.7).
\end{prf}

In the special case in which $X = \Spec \bbC$ there is the following
comparison result
\begin{prp}
  \label{prp:compcyc}
  Let $Y$ be a quasi-projective variety with projective closure
  $\overline{Y}$ and closed complement $Y^{\oo}$. Then the natural map
  $M_{\bt}(\Spec \bbC, \overline{Y}/Y^{\oo}) \to Z_{\bt}(Y,k)$ is a
  weak homotopy equivalence.
\end{prp}
\begin{prf}
  This follows from \cite{FriGab:Intersection} proposition 1.3 and
  theorem 1.4, which enable us to commute, up to homotopy equivalence,
  the singular chains functor with the group completion and the
  quotient in the definitions of $M_{\bt}(X,\overline{Y}/Y^{\oo}, k)$
  and $Z_{\bt}(Y,k)$.
\end{prf}

\begin{rmk}
  In particular, for a quasi-projective $Y$ with closed subvariety
  $Y_0$ and open complement $U$, theorem \ref{thm:loc} and proposition
  \ref{prp:compcyc} together with the long exact sequence of homotopy
  groups associated to a fibration gives the long exact sequence
  \begin{equation*}
    \xymatrix{\cdots \ar[r] &L_kH_n(Y_0) \ar[r]& L_kH_n(Y) \ar[r]& L_kH_n(U) \ar[r]&
      L_kH_{n-1}(Y_0) \ar[r]& \cdots}
  \end{equation*}
  This is the localization theorem for Lawson homology (proposition
  4.8 in \cite{Lim:Qproj}).
\end{rmk}

\begin{dfn}
  \label{dfn:susp}
  Let $Y$ be a projective variety with a given projective embedding $Y
  \subset \bbP^N$. Let $\bbP^{N} \subset \bbP^{N+1}$ be the embedding
  in the first $N$ coordinates, and $p^{\oo} = [0,\ldots,0,1] \in
  \bbP^{N+1}$ ($p^{\oo}$ is a choic of a point not in $\bbP^{N}$). The
  \emph{Lawson suspension} of $Y$, denoted by $\Lsu Y$ is the
  projective variety in $\bbP^{N+1}$ of the points contained in a line
  through $p^{\oo}$ that meets $Y$.
\end{dfn}

\begin{rmk}
  The Lawson suspension induces a continuous algebraic map at the
  level of Chow varieties
  \begin{equation*}
    \Lsu \cl C(Y,k) \to C(\Lsu Y,k+1),
  \end{equation*}
  such that 
  \begin{equation*}
    \Lsu(n_1 [W_1]+ \cdots + n_r[W_r]) =  n_1 [\Lsu W_1] + \cdots +
    n_r[\Lsu W_r].
  \end{equation*}
\end{rmk}

This leads to a very important theorem in the theory.
\begin{thm}[Lawson suspension]
  \label{thm:morsusp}
  The Lawson suspension of cycles induces a map of morphic spaces
  \begin{equation}
    \xymatrix{M_{\bt}(X,\overline{Y}/Y^{\oo},k) \ar[r] & M_{\bt}(X,\Lsu
    \overline{Y}/\Lsu Y^{\oo}, k+1)}
  \end{equation}
  which is a homotopy equivalence.
\end{thm}
\begin{prf}
  See \cite{FriLaw:Cocycles} theorem 3.3 for the projective case and
  \cite{Fri:CocyclesQproj} proposition 3.7 for the quasi-projective
  version.
\end{prf}

\begin{cor}
  \label{cor:LawHI}
  The projection $p\cl X\tms \bbA^1 \to X$ induces, through pull-back
  of cycles, a homotopy equivalence
  \begin{equation}
    Z_{\bt}(X, k) \weq Z_{\bt}(X\tms \bbA^1, k+1).
  \end{equation}
\end{cor}
\begin{prf}
  Using that $\Lsu X \setminus \set{p_{\oo}} = X\tms \bbA^1$, this
  becomes a consequence of the localization theorem \ref{thm:loc}, the
  suspension theorem \ref{thm:morsusp} and the homotopy equivalence
  \ref{prp:compcyc} between the space of cycles $Z_{\bt}(Y)$ and
  $M_{\bt}(\Spec \bbC, \overline{Y}/Y^{\oo}, k)$.
\end{prf}
\subsection{Morphic cohomology complexes}
\label{sec:morcpx}

As said in the begining of this section, to construct the morphic
cohomology we will need an algebraic analogue of $\Map(X,K(\bbZ,
2q))$. Since $Z(\bbA^q,0)$ is an algebraic model of a $K(\bbZ,
2q)$, it would makes sense to use the spaces $M_{\bt}(X,\bbA^t,k)$. We
have defined those spaces for a pair of projective varieties
$(\overline{Y}, Y^{\oo})$ not for a quasi-projective $Y$ alone, but
$\bbA^t$ has a natural choice of projective closure,
$\overline{Y}=\bbP^t$ for which the closed complement is
\begin{equation*}
  Y^{\oo} = \set{[x_0:\cdots:x_t] \in \bbP^t \mid x_0 = 0} =
  \bbP^{t-1}.
\end{equation*}
So a good candidate spaces computing the morphic cohomology are
$M_{\bt}(X, \bbP^t/\bbP^{t-1}, k)$. There are two indices involved in
this construction, $t$ and $k$ apart from the one appearing when
taking homotopy groups, but thanks to the suspension theorem
\ref{thm:morsusp}, one index amongst $t$ and $k$ is redundant as
\begin{equation*}
  M_{\bt}(X, \bbP^t/\bbP^{t-1}, k) \weq M_{\bt}(X, \bbP^{t+1}/\bbP^{t}, k+1)
\end{equation*}
is a homotopy equivalence. What really matters is the difference
$q=t-k$. So we get rid of the useless index in a way that ensures good
functorial properties in situations where the suspension is involved.
\begin{dfn}
  \label{dfn:mor}
  The \emph{morphic space} of level $q$ associated to a
  quasi-projective algebraic variety $X$ is the simplicial abelian
  group $M_{\bt}(X,q)$ given as the colimit
  \begin{equation}
    M_{\bt}(X,q) = \colim_{k} M_{\bt}(X, \bbP^{q+k}/\bbP^{q+k-1}, k).
    \label{eq:morspc}
  \end{equation}
  where the maps 
  \begin{equation*}
    M_{\bt}(X, \bbP^{q+k}/\bbP^{q+k-1}, k) \weq M_{\bt}(X, 
    \bbP^{q+k+1}/\bbP^{q+k}, k+1)
  \end{equation*}
  are the homotopy equivalences in the suspension theorem.

  The \emph{morphic cohomology groups} of $X$ are the groups
  $L^qH^n(X)$ given by
  \begin{equation}
    L^qH^n(X) = \pi_{2q-n}M_{\bt}(X, q).
    \label{eq:dfnmor}
  \end{equation}
\end{dfn}

\begin{rmk}
  All the maps $M_{\bt}(X, \bbP^{q+k}/\bbP^{q+k-1}, k) \to
  M_{\bt}(X,q)$ are homotopy equivalences. In particular,
  \begin{equation}
    L^qH^n(X) \simeq \pi_{2q-n} M_{\bt}(X, \bbP^{q+k}/\bbP^{q+k-1}, k)
  \end{equation}
  for any $k \geq 0$.
\end{rmk}

Recall that, for an abelian category $\caA$, the Dold-Kan
correspondence establishes an equivalence
\begin{equation}
  N\cl s\caA \to \ctCh_+(\caA)
\end{equation}
between the categories $s\caA$ of simplicial objects in $\caA$ and the
category $\ctCh_+(\caA)$ of positively graded chain complexes in
$\caA$. For every simplicial object $A_{\bt}$ in $s\caA$, the
associated normalized chain complex $N_*(A)$ is such that
\begin{equation}
  \pi_n A_{\bt} \simeq H_n N_*(A).
\end{equation}

Applying the Dold-Kan functor $N$ to the morphic cohomology spaces
$M_{\bt}(X,q)$ we obtain a chain complex of abelian groups. It is
convenient, in analogy with the motivic world, to reindex those
morphic complexes as follows:
\begin{equation}
  M^n(X, q) = M_{2q-n}(X, q),
  \label{eq:morcpx}
\end{equation}
obtaining a cochain complex of abelian groups $M^*(X, q)$, which is
zero for $n > 2q$, so it is bounded above, but is unbounded below.

There is a couple of important properties we will use.
\begin{thm}
  \label{thm:morHI}
  Let $X$ be an algebraic variety and $p\cl E \to X$ be a vector
  bundle over $X$. Then the induced map
  \begin{equation}
    M_{\bt}(X, q) \to M_{\bt}(E,q)
  \end{equation}
  is a weak homotopy equivalence and, in particular, the respective
  morphic cohomology groups are isomorphic.
\end{thm}
\begin{prf}
  Proposition 3.5 in \cite{Fri:CocyclesQproj} tells us that the map
  $M_{\bt}(X, Y,q) \to M_{\bt}(E, Y,q)$ between bivariant spaces is a
  homotopy equivalence for $Y$ projective. Then a five lemma argument
  applied to the pair $\bbP^q/\bbP^{q-1}$ proves what we want.
\end{prf}

\begin{thm}[Duality]
  \label{thm:morduality}
  Let $X$ be a quasi-projective variety of dimension $d$. There
  is a natural map
  \begin{equation}
    \Gamma\cl M_{\bt}(X, q) \to Z_{\bt}(X\tms\bbA^q,d)
    \label{eq:graphing}
  \end{equation}
  called the \emph{graphing construction}.

  Moreover, when $X$ is smooth, $\Gamma$ becomes a homotopy
  equivalence. In such case gives isomorphisms
  \begin{equation}
    L^qH^n(X) \weq L_{d-q}H_{2d-n}(X).
    \label{eq:duality}
  \end{equation}
\end{thm}
\begin{prf}
  The map $\Gamma$ is constructed and proved to be a homotopy
  equivalence for smooth varieties in \cite{FriLaw:Duality} theorem
  3.3 when $X$ is projective and extended to the quasi-projective case
  in \cite{Fri:CocyclesQproj} theorem 5.2. The isomorphism
  \eqref{eq:duality} comes from the combination of the homotopy
  equivalence $\Gamma$ and the homotopy invariance for cycle spaces
  \ref{cor:LawHI}.
\end{prf}

\begin{cor}
  \label{cor:smmorMV}
  Let $X$ be smooth and $\set{U,V}$ an open covering of $X$. Then the
  diagram
  \begin{equation}
    \xymatrix{M^*(X,q) \ar[r] \ar[d] & M^*(U,q) \ar[d]\\
      M^*(V,q) \ar[r] & M^*(U\cap V, q)}
    \label{eq:smmorMVsq}
  \end{equation}
  gives rise to a Mayer-Vietoris long exact sequence for morphic
  cohomology
  \begin{equation}
    \xymatrix{\cdots \ar[r] & L^qH^n(X) \ar[r] &
      L^qH^n(U)\oplus L^qH^n(V) \ar[r] &
      L^qH^n(U\cap V) \ar[r] & \cdots}
    \label{eq:smmorMV}
  \end{equation}
\end{cor}
\begin{prf}
  As $X$ is smooth, The duality theorem implies that the square
  \eqref{eq:smmorMVsq} is quasi-isomorphic to the diagram obtained
  applying the Dold-Kan functor to
  \begin{equation*}
    \xymatrix{Z_{\bt}(X\tms\bbA^q,d) \ar[r] \ar[d] & Z_{\bt}(U\tms\bbA^q,d) \ar[d]\\
      Z_{\bt}(V\tms\bbA^q,d) \ar[r] & Z_{\bt}(U\cap V\tms\bbA^q, d)}
  \end{equation*}
  That last square is homotopy cartesian, because by the localization
  theorem \ref{thm:loc} the homotopy fibers of the horizontal maps are
  $Z_{\bt}((X\setminus U)\tms \bbA^q,d)$ and $Z_{\bt}((V\setminus
  U)\tms \bbA^q, d)$, and those spaces are homotopy equivalent (in
  fact isomorphic) as $X\setminus U = V\setminus U$.
\end{prf}

\subsection{A modification of the morphic cohomology}
\label{sec:mormod}
We are interested in using homotopy invariance and Mayer-Vietoris
properties for morphic cohomology. Those properties hold for smooth
varieties thanks to duality (theorems \ref{thm:morHI} and
\ref{cor:smmorMV}) but unfortunately are not known (at least to the
author) to hold for singular varieties. So we modify the definition of
morphic cohomology to force those two properties. This is not new. In
\cite{Fri:BlochOgus} Friedlander defines a modified morphic
cohomology, called topological cycle cohomology, which satisfies a
Mayer-Vietoris propery for Zariski open covers. In fact, the point of
\cite{Fri:BlochOgus} is to prove that this theory together with a
homological companion satisfies the Bloch-Ogus axioms.

First of all, note that the contravariant functoriality of the
complexes $M^*(X, q)$ with respect to $X$ makes them a cochain
complexes of presheaves of abelian groups on the category
$\ctqProj_\bbC$.

Let $\Delta_{\bt}$ be the standard cosimplicial scheme
\begin{equation*}
  \Delta_n = \set{(x_0,\ldots,x_n) \in \bbA^{n+1} \mid x_0+\cdots +
    x_n = 1}.
\end{equation*}
\begin{dfn}
  \label{dfn:morshv}
  We define the \emph{morphic complexes of sheaves} $\caM_{\suzar}^*(q)$ as the
  Zariski sheafification of the presheaf of complexes $M^*(\hy, q)$,
  that is
  \begin{equation}
    \caM_{\suzar}^n(q) =  \mrm{sh}_{\suzar}M^*(\hy,q).
    \label{eq:morshv}
  \end{equation}

  Additionally, we define the \emph{homotopy invariant morphic
  complexes} as the following total complex
\begin{equation}
  \caM_{\suhi}^n(q) = \Tot^n \caM_{\suzar}^*(\hy\tms \Delta_\bt,q).
  \label{eq:morshvHI}
\end{equation} 
\end{dfn}

\begin{rmk}
  \label{rmk:resolutions}
  The complexes $\caM^*(q)$ are unbounded below and bounded above by
  $2q$. This poses some homological algebra troubles as they are
  bounded on the wrong side. However, due to a result of Spaltenstein
  \cite{Spa:UnbCpx} one can still have resolutions $\caM^*(q) \to
  \caI^*$ playing the role of injective resolutions. Those are called
  K-injective in \cite{Spa:UnbCpx}. Later, we will use cohomological
  finiteness arguments to prove convergence of the spectral sequences
  we encounter.
\end{rmk}

We will think of the complexes $\caM_{\suzar}^*(q)$ and
$\caM_{\suhi}^*(q)$ as objects in the derived category
$\mbf{D}^-\ctSh(\ctqProj_{\bbC})$ of abelian sheaves on the Zariski
site $\ctqProj_{\bbC}$.

For any variety $X$, and a complex of sheaves $\caF^* \in
\mbf{D}^-\ctSh(\ctqProj_{\bbC})$, we will denote by $\caF^*|_X$ the
restriction of $\caF^*$ to the small zariski site of $X$,
i.e. $\caF^*|_X$ is an element in the category $\mbf{D}^-\ctSh(X)$,
the derived category of abelian sheaves on $X$.

This sheaf-theoretic interpretation leads to natural reformulations of
the morphic cohomology groups.
\begin{dfn}
  \label{dfn:dmor}
  The \emph{topological cycle cohomology groups} are the
  hypercohomology groups
  \begin{equation}
    L^qH_{\suzar}^n(X) = \bbH^n(X, \caM_{\suzar}(q)|_X).
    \label{eq:zarmor}
  \end{equation}

  The \emph{homotopy invariant morphic cohomology} are the
  hypercohomology groups
  \begin{equation}
    L^qH_{\suhi}^n(X) = \bbH^n(X, \caM_{\suhi}(q)|_X).
    \label{eq:himor}
  \end{equation}
\end{dfn}
\begin{rmk}
  The topological cycle cohomology groups $L^qH_{\suzar}^n(X)$ are the
  ones defined by Friedlander in \cite{Fri:BlochOgus}. Although the
  original morphic cohomology groups are homotopy invariant by
  \ref{thm:morHI}, we are not able to prove this homotopy invariance
  for $L^qH_{\suzar}^n$. This is the reason for introducing the theory
  $L^qH_{\suhi}^n$ which exhibits both properies, homotopy invariance
  and Mayer-Vietoris.
\end{rmk}
\begin{rmk}
  The hypercohomologies in \ref{dfn:dmor} can be rephrased as an
  ext-group on the derived category $\mbf{D}^-\ctSh(X)$ of sheaves of
  abelian groups on $X$, for example
  \begin{equation}
    L^qH_{\suhi}^n(X) = \Ext^n_{\ctSh(X)}(\bbZ_X, \caM_{\suhi}^*(q)|_X),
  \end{equation}
  where $\bbZ_X$ denotes the constant sheaf on $X$.

  Moreover, let $\bbZ h_X$ denote the free abelian group generated by
  the sheaf on $\ctqProj_\bbC$ represented by $X$, i.e. $\bbZ h_X(U) =
  \bbZ \Hom(U,X)$. Then there is an isomorphism
  \begin{equation}
    L^qH_{\suhi}^n(X) = \Ext^n_{\ctSh(\ctqProj_{\bbC})}(\bbZ h_X, \caM_{\suhi}^*(q)),
  \end{equation}
  where now the $\Ext$ is taken in the category of sheaves on the big
  site $\ctqProj_{\bbC}$.

  We will use either interpretation, as an $\Ext$ in the big site or
  in the small site.
\end{rmk}

This new versions of morphic cohomology are related to the old one as
follows. There are comparison morphisms
\begin{equation}
  \xymatrix{M^*(\hy, q) \ar[r]^{a}& \caM_{\suzar}^*(q) \ar[r]^b & \caM_{\suhi}^*(q)},
\end{equation}
the first is the inclusion of a presheaf in its associated sheaf, and
the second is the inclusion into the summand of the total complex
corresponding to the algebraic simplex $\Delta_0$.

As defined in \ref{dfn:mor}, $L^qH^n(X)$ is the homology of the
complex of global sections of the cochain complex of presheaves
$M^*(\hy, q)$, while the other two flavours of morphic cohomology are
the hypercohomology of the respective sheaves.

Let's take then K-injective resolutions $\caM_{\suzar}^*(q) \to
\caI_{\suzar}^*$ and $\caM_{\suhi}^*(q) \to \caI_{\suhi}^*$. There are
induced comparison maps
\begin{equation*}
  \xymatrix{M^*(\hy, q) \ar[r]^{a}& \caI_{\suzar}^* \ar[r]^b & \caI_{\suhi}^*},
\end{equation*}
which, when taking cohomology of the global sections induce comparison
maps
\begin{equation}
  \xymatrix{L^qH^n(X) \ar[r]^{a}& L^qH^n_{\suzar}(X) \ar[r]^b &
    L^qH^n_{\suhi}(X)}.
  \label{eq:morcomp}
\end{equation}

\begin{prp}
  \label{prp:smzardesc}
  Assume $X$ is a smooth quasi-projective variety. Then the comparison
  morphisms \eqref{eq:morcomp} are isomorphisms.
\end{prp}
\begin{prf}
  From \ref{cor:smmorMV} we know that the presheaf $M_{\bt}(\hy, q)$
  restricted to smooth varieties satisfies the Mayer-Vietoris propery,
  that is, the square
  \begin{equation*}
    \xymatrix{M^*(X,q) \ar[r] \ar[d] & M^*(U,q) \ar[d]\\
              M^*(V,q) \ar[r] & M^*(U\cap V, q)}
  \end{equation*}
  is homotopy cartesian. Take a K-injective resolution
  $\caM_{\suzar}^*(q)\to \caI_{\suzar}^*$. Using the Brown-Gersten
  theorem (\cite{Voe:CD} lemma 3.5) we conclude that the presheaf
  $M^*(\hy, q)$ is globally weakly equivalent to the K-injective
  resolution $\caI^*$, so
  \begin{equation*}
    H^n\Gamma(X, M^*(\hy, q)) \simeq H^n\Gamma(X, \caI_{\suzar}^*),
  \end{equation*}
  and this settles the first isomorphism $a$.

  As for $b$, by theorem \ref{thm:morHI} we know that $L^qH^n(\hy)$ is
  a homotopy invariant functor, and by the previous isomorphism
  coincides with $L^qH^n_{\suzar}(\hy)$ on smooth varieties, so the
  last one is also homotopy invariant on smooth varieties. Then the
  spectral sequence associated to the double complex defining
  $L^qH^n_{\suhi}$ degenerates on the second page and this gives the
  isomorphism $b$.
\end{prf}

Finally, those are the properties we will use of the modified theory
$L^qH^n_{\suhi}$.
\begin{thm}
  \label{thm:morhiMV}
  Let $U,V$ be an open cover of a complex variety $X$. Then there is a
  Mayer-Vietoris long exact sequence for the groups $L^qH^n_{\suhi}$
\begin{equation}
\xymatrix{\cdots \ar[r] & L^qH_{\suhi}^n(X) \ar[r] &
  L^qH_{\suhi}^n(U)\oplus L^qH_{\suhi}^n(V) \ar[r] &
  L^qH_{\suhi}^n(U\cap V) \ar[r] & \cdots}
\end{equation}
\end{thm}
\begin{prf}
  The open cover $\set{U,V}$ gives an exact triangle in
  $\mbf{D}^-\ctSh(X)$
  \begin{equation}
  \xymatrix{\bbZ_{U\cap V} \ar[r]& \bbZ_{U} \oplus \bbZ_{V} \ar[r] &
    \bbZ_X \ar[r]& \bbZ_{U\cap V}[1] &}
  \end{equation}
  which, when passed through the functor $\Ext^n(\hy,
  \caM_{\suhi}^*(q))$, gives the desired long exact sequence.
\end{prf}

\begin{thm}
  \label{thm:morhiHI}
  The groups $L^qH^n(X)$ are homotopy invariant, that is the
  projection $p\cl X\tms \bbA^1 \to X$ induces isomorphisms
  \begin{equation}
    L^qH^n(X) \weq L^qH^n(X\tms \bbA^1).
  \end{equation}
\end{thm}
\begin{prf}
  This is a standard argument. See for example corollary 2.19 in \cite{MVW:motivic}
\end{prf}

\subsection{Some remarks about notation}
\begin{rmk}
  In the remainder of this paper we will deal only with the version
  $L^qH^n_{\suhi}$ of morphic cohomology, as we are interested in
  using a version of a Mayer-Vietoris spectral sequence and homotopy
  invariance. We will drop the subindex ``$\suhi$'' from the notation for
  readability.
\end{rmk}

It will be useful to deal with the morphic complexes all at once, so
we define the complex of sheaves
\begin{equation*}
  \caM^* = \bigoplus_{q \geq 0} \caM^*(q).
\end{equation*}
This is a bigraded object, having the algebraic degree $q$ and the
cohomological degree $n$ with a differential for the $n$-grading. We
will denote by $\caH^n\caM$ its cohomology sheaves.

Analogously, we will use the notation $LH^n$ for the graded group
\begin{equation*}
  LH^n(X) = \bigoplus_{q\geq 0} L^qH^n(X) = \bbH^n(X, \caM^*),
\end{equation*} 
and $LH$ for the bigraded groups
\begin{equation*}
  LH(X) = \bigoplus_{n \geq 0} LH^n(X).
\end{equation*}

\subsection{Cup product and Kunneth homomorphism}
\label{sec:cup}
Let $X,X'$ be schemes. Following \cite{FriWa:Function} proposition
3.2, the projections $\pi_X \cl X\tms X' \to X$ and $\pi_{X'}\cl X
\tms X' \to X'$ induce exterior products on morphic complexes
\begin{equation}
  M^*(X,q)\tsr M^*(X',q') \to M^*(X \tms X',q+q').
  \label{eq:extprodcpx}
\end{equation}

Composing this exterior product with the pull-back through the
diagonal embedding $\Delta\cl X \to X\tms X$ we get a cup product
\begin{equation}
  M^*(X,q)\tsr M^*(X,q') \to M^*(X \tms X,q+q') \to M^*(X,q+q'),
\end{equation}
which induces a cup product
\begin{equation}
  \caM^*(q)\tsr \caM^*(q') \to \caM^*(q+q').
  \label{eq:cupsh}
\end{equation}
on the sheaves $\caM^*$.

\begin{thm}
  \label{thm:cup}
  The cup porduct \eqref{eq:cupsh} induces an associative product on
  morphic cohomology
  \begin{equation}
    L^qH^n(X)\tsr L^{q'}H^{n'}(X) \to L^{q+q'}H^{n+n'}(X).
  \end{equation}
  
  which is graded commutative with respect to the cohomological
  grading, that is, for $a\in LH^n(X)$, $b \in LH^{n'}(X)$ the
  commutativity relation
  \begin{equation*}
    a \cup b = (-1)^{nn'} b \cup a
  \end{equation*}
  holds.

  Moreover, the cup product is functorial in the sense that for $f\cl
  X \to Y$ a morphism of algebraic varieties, the induced map on
  morphic cohomology $f^*\cl LH(Y) \to LH(X)$ is a ring homomorphism.
\end{thm}
\begin{prf}
  See \cite{FriLaw:Cocycles} corollary 6.2.
\end{prf}

Let $\Mpt = \Gamma(\Spec \bbC, \caM^*)$ be the sections of the morphic
complex on a point, and let $\LHpt$ be its hypercohomology, which in
this case coincides with the cohomology ring of the global sections
$\Mpt$.
\begin{cor}
  \label{cor:module}
  The morphic cohomology of a complex quasi-projective variety $X$ has
  a canonical structure of $\LHpt$-module.
\end{cor}
\begin{prf}
  Applying the second part of theorem \ref{thm:cup} to the structure
  map $X \to \Spec \bbC$ we get a ring homomorphism $\LHpt \to LH(X)$,
  making $LH(X)$ into a $LHpt$-module.
\end{prf}

Now, the exterior product \eqref{eq:extprodcpx} also induces an
exterior product on morphic cohomology
\begin{equation}
  L^qH^n(X)\tsr L^{q'}H^{n'}(X') \to L^{q+q'}H^{n+n'}(X\tms X').
  \label{eq:extprod}
\end{equation}
As the structure of $LHpt$-module in $LH(X)$ is functorial, the action
of $\LHpt$ on either factor on the left of \eqref{eq:extprod}, gives
the same action on the right, so the exterior product map
\eqref{eq:extprod} factors through
\begin{equation}
  LH(X)\tsr_{LHpt}LH(Y) \to LH(X\tms Y).
  \label{eq:kunnethmap}
\end{equation}
In rather special circumstances this Kunneth homomorphism
\eqref{eq:kunnethmap} happens to be an isomorphism. We are interested
in a very special case of this Kunneth isomorphism, which we now
prove.
\begin{prp}
  \label{prp:kunneth}
  Let $X$ be a smooth quasi-projective variety. The Künneth
  homomorphism
  \begin{equation}
    LH(X)\tsr_{LHpt}LH(\bbG_m) \to LH(X\tms \bbG_m).
  \end{equation}
  is an isomorphism.
\end{prp}
\begin{prf}
  Let $i\cl pt \to \bbA^1$ be the inclusion of a point and $j\cl
  \bbG_m\to \bbA^1$ its open complement. We have the following
  commutative diagram of long exact sequences

  \begin{equation*}
    \xymatrix{\cdots LH(X)\tsr_{LHpt} LH(pt) \ar[r]^{id\tms i_!} \ar[d]&
      LH(X)\tsr_{\LHpt} LH(\bbA^1) \ar[r]^{id\tms j^*} \ar[d]&
      LH(X)\tsr_{\LHpt} LH(\bbG_m) \ar[d] \cdots\\
      \cdots LH(X\tms pt) \ar[r]^{(id\tms i)_!} & LH(X\tms \bbA^1)
      \ar[r]^{(id\tms j)^*} & LH(X\tms
      \bbG_m) \cdots}
  \end{equation*}
  The vertical maps are the Kunneth morphisms, and $i_!$ is the Gysin
  map defined by duality (theorem \ref{thm:morduality}) as $i_! =
  \Gamma_*^{-1}i_*\Gamma_*$. The exactness of the rows come, by
  duality, from the localization theorem \ref{thm:loc} and the long
  exact sequence of homotopy groups of a fibration.

  Now a standard application of the five lemma together with the
  homotopy invariance \ref{thm:morHI} proves the desired isomorphism.
\end{prf}

\begin{rmk}
  I became aware of a construction of a Kunneth spectral sequence in
  the case in which $X$ or $Y$ is a linear variety in a private
  communication with Mircea Voineagu \cite{Voi:Priv}. An analogous
  construction is done in \cite{Jos:KthLinear} for higher Chow groups and
  K-theory.
\end{rmk} 

Finally, some basic computations of morphic cohomology rings we will need.
\begin{prp}
  \label{prp:comp}
  \begin{enumerate}
  \item For $k \geq 0$,
    \begin{equation*}
      LH(\bbA^k) \simeq \bbZ[s],
    \end{equation*}
    where $s$ is a free generator of bidegree $(1,0)$ (degree 1 with
    respect to the $q$-grading).

  \item For $k \geq 0$,
    \begin{equation*}
      LH(\bbP^k) \simeq \bbZ[s,h]/(h^{k+1})
    \end{equation*}
    where $s$ has bidegree $(1,0)$ and $h$ has bidegree $(1,2)$.

  \item The morphic cohomology of the multiplicative group is given by
    \begin{equation*}
      LH(\bbG_m)\simeq \bbZ[s,e]/(e^2),
    \end{equation*}
    where $s$ is a generator of bidegree $(1,0)$ and $e$ is a
    generator of bidegree $(1,1)$.
  \end{enumerate}
\end{prp}
\begin{prf}
  1) and 2) follow from duality and the computations of Lawson
  homology using the suspension theorem \ref{thm:morsusp}.

  3) Take the open cover of $\bbP^1$ by two affines. Then we have the
  following piece of Mayer-Vietoris sequence
\begin{equation*}
  \xymatrix{LH^n(\bbP^1) \ar[r] & LH^n(\bbA^1)^{\oplus 2} \ar[r] &
    LH^n(\bbG_m) \ar[r]& LH^{n+1}(\bbP^1) \ar[r] & LH^{n+1}(\bbA^1)^{\oplus2}}
\end{equation*}
which, using 1) and 2) for the computations of $\bbP^1$ and $\bbA^1$
gives the result.
\end{prf}

\begin{rmk}
  Note that, in particular, $\LHpt \simeq \bbZ[s]$. As a consequence
  of \ref{cor:module}, the morphic cohomology ring $LH(X)$ is, in
  fact, a $\bbZ[s]$-module. The action by $s$ on $LH(X)$ corresponds
  to the $s$-maps defined in \cite{FriMaz:Operations}.
\end{rmk}

\section{Toric varieties}
\label{sec:toric}
First we set the notation. Let $N \simeq \bbZ^n$ be a free
$\bbZ$-module of rank $n$, and $M$ its dual $\bbZ$-module. We will
denote by $N_{\bbR} = N\tsr\bbR$ and $M_{\bbR}=M\tsr\bbR$. In this
way, $N$ and $M$ are to be thought as lattices on $N_{\bbR}$ and
$M_{\bbR}$. Moreover, there is the duality pairing $\angbk{u,v}$ for
$u \in N_{\bbR}$ and $v\in M_{\bbR}$.

\subsection{Cones and fans}
\begin{dfn}
  \label{dfn:cone}
  A \emph{rational polyhedral cone} in $N$ is a set $\sigma \subset
  N_{\bbR}$ generated by a finite number of integral vectors
  $v_1,\ldots,v_k \in N$ in the following way:
  \begin{equation*}
    \sigma = \set{\lambda_1v_1 + \cdots + \lambda_k v_k \in N_{\bbR}
      \mid \lambda_i \in \bbR_{\geq 0}}.
  \end{equation*}

  The \emph{dimension} of a cone $\sigma$ is the dimension of the
  vector space $\bbR\sigma$ generated by the cone.

  If a rational polyhedral cone $\sigma$ does not contain a vector
  space of positive dimension is said to be \emph{strictly convex}.
\end{dfn}

\begin{dfn}
  \label{dfn:dual}
  Let $\sigma$ be a rational polyhedral cone in $N$.
  \begin{enumerate}
  \item The \emph{dual cone} of $\sigma$ is a rational polyhedral cone
    in $M$ given by
    \begin{equation*}
      \sigma^{\vee} = \set{v \in M_{\bbR} \mid \angbk{u,v} \geq 0
        \text{ for all $u \in \sigma$}}.
    \end{equation*}
  \item The \emph{orhogonal cone} of $\sigma$ is a rational polyhedral
    cone in $M$ given by
    \begin{equation*}
      \sigma^{\bot} = \set{v \in M_{\bbR} \mid \angbk{u,v} = 0 \text{
          for all $u \in \sigma$}}.
    \end{equation*}
  \end{enumerate}
\end{dfn}

\begin{rmk}
  The orthogonal cone $\sigma^{\bot}$ is in fact a vector space, as if
  $\angbk{u,v}=0$ then $\angbk{u,\lambda v} = 0$ for all $\lambda \in
  \bbR$.
\end{rmk}

\begin{dfn}
  \label{dfn:face}
  A cone $\tau$ is a \emph{face} of $\sigma$ if there exists a $u \in
  M_{\bbR}$ such that
  \begin{equation*}
    \tau = \sigma \cap u^{\bot}.
  \end{equation*}
\end{dfn}

\begin{dfn}
  \label{dfn:fan}
  A \emph{fan} $\Delta$ is a set of strongly convex rational
  polyhedral cones in $N$ such that
  \begin{enumerate}
  \item Every face of a cone in $\Delta$ also belongs to $\Delta$.
  \item The intersection of two cones in $\Delta$ is a face of both
    intersecting cones.
  \end{enumerate}

  The notation $\Delta^{(k)}$ will mean the set of all cones of
  codimension $k$ in $\Delta$.
\end{dfn}
\begin{rmk}
  A fan has a partial order given by the inclusion of faces. We will
  use the notation $\tau \leq \sigma$ to say that $\tau$ is a face of
  $\sigma$.
\end{rmk}

\begin{dfn}
  \label{dfn:ori}
  An \emph{orientation of a cone} $\sigma$ is an orientation of the
  vector spaces $\bbR \sigma$.

  An \emph{orientation of a fan} $\Delta$ will be a choice of an
  orientation for every cone in $\Delta$.
\end{dfn}

We will always use fans with a fixed orientation.

\begin{rmk}
  Let $\tau \leq \sigma$ be a face of codimension 1 in $\sigma$. This
  means that there exists a $u \in \sigma^{\vee}$ such that $\tau =
  \sigma \cap u^{\bot}$. Then we have an isomorphism
  \begin{equation}
    \bbR \sigma \simeq \bbR u \oplus \bbR\tau.
    \label{mth:oriso}
  \end{equation}
  which allows us to transfer the orientation of $\sigma$ to $\tau$ as
  follows: the orientation induced on $\tau$ by $\sigma$ is the one
  compatible with the isomorphism \fref{mth:oriso} and taking the
  orientation in $\bbR u$ given by the vector $u$.
\end{rmk}

\subsection{Definition of a toric variety}
Let us fix a field $k$. For now it will be an aribrary field, but in
the next sections we will set $k=\bbC$.

We describe the toric variety $X(\Delta)$ associated to a fan
$\Delta$. The scheme $X(\Delta)$ is constructed locally, one affine
piece for every cone in $\Delta$, and then glueing them, according to
the combinatorics of the fan, to build the scheme $X(\Delta)$.

For every cone $\sigma \in \Delta$ we define an affine scheme
\begin{equation*}
  X_{\sigma} = \Spec k[\sigma^{\vee}\cap M],
\end{equation*}
where $\sigma^{\vee}\cap M$ is seen as a monoid and
$k[\sigma^{\vee}\cap M]$ is its monoid algebra. An inclusion of cones
$\tau \leq \sigma$ induces an inclusion on dual cones $\sigma^{\vee}
\leq \tau^{\vee}$ and this induces a morphism on the corresponding
monoid algebras, obtaining a morphism of schemes $i_{\tau,\sigma}\cl
X_{\tau} \to X_{\sigma}$.

\begin{prp}
  \label{prp:toric}
  Let $\Delta$ be a fan on a lattice $N$ of dimension $n$.
  \begin{enumerate}
  \item For every face inclusion $\tau \leq \sigma$, the induced
    morphism $i_{\tau,\sigma}\cl X_{\tau}\to X_{\sigma}$ is an open
    embedding.

  \item There is a (unique) $k$-scheme $X(\Delta)$ together with a
    family of open embeddings $i_{\sigma}\cl X_{\sigma} \to X(\Delta)$
    for every cone in $\Delta$ such that $X_{\sigma}$ is an open
    covering of $X(\Delta)$ and for every pair of cones $\tau \leq
    \sigma$, the following diagram commutes
    \begin{equation*}
      \xymatrix{X_{\tau} \ar[r]^{i_{\tau,\sigma}} \ar[rd]_{i_{\tau}} &
      X_{\sigma} \ar[d]^{i_{\sigma}}\\
      & X(\Delta)}
    \end{equation*}

  \item The scheme $X(\Delta)$ is reduced, normal and
    Noetherian. Moreover, if the cones in $\Delta$ cover the whole
    $N_{\bbR}$ the variety $X(\Delta)$ is complete.

  \item For $\tau,\sigma \in \Delta$, the face inclusions of $\tau
    \cap \sigma$ induce an isomorphism
    \begin{equation*}
      X_{\tau\cap\sigma} \simeq X_{\tau}\times_{X(\Delta)} X_{\sigma}.
    \end{equation*}
  \end{enumerate}
\end{prp}
\begin{prf}
  See \cite{Ful:Toric}, sections 1.4 and 2.1.
\end{prf}

\begin{dfn}
  \label{dfn:toric}
  A \emph{toric variety} for a fan $\Delta$ is an algebraic variety
  $X(\Delta)$ associated $\Delta$ as in the previous theorem.
\end{dfn}

\begin{exm}
  Let $N$ be a lattice of rank $n$ and $\Delta$ the fan containing
  only the cone $\set{0}$. Then $X(\Delta)$ is an algebraic torus of
  dimension $n$, that is there is an isomorphism
  \begin{equation*}
    X(\Delta) = \Spec k[M] \simeq \bbG_m^n.
  \end{equation*}
  This algebraic torus is, in fact, a commutative group scheme, and
  the group law is described by the comultiplication on the coordinate
  ring $k[M] \to k[M] \tsr k[M]$ given by $m \mapsto m\tsr m$ for
  every $m \in M$ and extended by linearity.
\end{exm}

\begin{rmk}
  Every toric variety $X(\Delta)$ contains an algebraic torus as a
  dense open set, corresponding to the affine piece
  \begin{equation*}
    X_{0} = \Spec k[M] \simeq \bbG_m^{n}
  \end{equation*}
  associated to the cone ${0} \in \Delta$.

  The group law on the torus $X_{0}$ can be extended to an action on
  the whole toric variety. Locally on an open piece $X_{\sigma} =
  \Spec k[\sigma^{\vee}\cap M]$ this action is given by the morphism
  of rings
  \begin{equation*}
    k[\sigma^{\vee}\cap M] \to k[M] \tsr k[\sigma^{\vee}\cap M]
  \end{equation*}
  such that $m \to m\tsr m$ for every $m \in \sigma^{\vee}\cap
  M$. Then all those local actions glue together forming an action
  \begin{equation*}
    X_0 \tms X(\Delta) \to X(\Delta).
  \end{equation*}
\end{rmk}

To every cone $\sigma \in \Delta$, we can associate an affine scheme
$T_{\sigma}$ besides the open affine $X_{\sigma}$ we have already
constructed,
\begin{equation*}
  T_{\sigma} = \Spec k[\sigma^{\bot}\cap M].
\end{equation*}

Observe that There is a homomorphism of rings
\begin{equation*}
  k[\sigma^{\vee}\cap M] \to k[\sigma^{\bot}\cap M]
\end{equation*}
sending every element of $\sigma^{\bot}$ to themselves and the
elements of $\sigma^{\vee}$ not in $\sigma^{\bot}$ to 0. This
homomorphism induces a morphism of schemes
\begin{equation*}
  T_{\sigma}=\Spec k[\sigma^{\bot}\cap M] \to X_{\sigma}=\Spec 
  k[\sigma^{\vee}\cap M],
\end{equation*}
and composing with the embedding $i_{\sigma}\cl X_{\sigma}\to
X(\Delta)$ we get morphisms $j_{\sigma}\cl T_{\sigma}\to X(\Delta)$.

On the other hand, we have retractions $r_{\sigma}\cl X_{\sigma} \to
T_{\sigma}$ defined by the inclusion $\sigma^{\bot} \to
\sigma^{\vee}$.

\begin{prp}
  \label{prp:orbits}
  For every cone $\sigma \in \Delta$, the morphism
  \begin{equation*}
    j_{\sigma}\cl T_{\sigma} \to X(\Delta),
  \end{equation*}
  is a locally closed embedding and corresponds with an orbit of the
  torus action.

  Moreover, the collection of all $T_{\sigma}$ for all the cones
  $\sigma \in \Delta$ give the orbit decomposition of the torus action
  on the toric variety $X(\Delta)$.
\end{prp}
\begin{prf}
  See \cite{Ful:Toric} section 3.1.
\end{prf}

\begin{rmk}
  Note that the orbits $T_{\sigma}$ are algebraic tori whose dimension
  is the codimension of $\sigma$, i.e.
  \begin{equation*}
    \dim T_{\sigma} = n-\dim\sigma.
  \end{equation*}
\end{rmk}

In short, associated to a toric variety we have an open cover given by
the affine pieces
\begin{equation*}
  X_{\sigma} = \Spec k[\sigma^{\vee}\cap M],
\end{equation*}
and an orbit decomposition as a disjoint union of tori
\begin{equation*}
  T_{\sigma} = \Spec k[\sigma^{\bot}\cap M].
\end{equation*}

There is a relation between them in terms of the combinatorics of the
fan as follows:
\begin{prp}
  \label{prp:desc}
  We have the following disjoint union decompositions,
  \begin{equation*}
    X_{\sigma} = \bigcup_{\tau \leq \sigma} T_{\tau},
  \end{equation*}
  \begin{equation*}
    \overline{T_{\tau}} = \bigcup_{\tau \leq \sigma} T_{\sigma}.
  \end{equation*}
\end{prp}
\begin{prf}
  See \cite{Ful:Toric} section 3.1.
\end{prf}

\begin{rmk}
  In particular, the orbit $T_{\sigma}$ is closed in $X_{\sigma}$, and
  is the lowest dimensional orbit contained in it.
\end{rmk}


\begin{prp}
  \label{prp:htpy}
  Let $X(\Delta)$ be a toric variety and $\sigma$ a cone in
  $\Delta$. There is a morphism
  \begin{equation*}
    h\cl X_{\sigma} \times \bbA^1_k \to X_{\sigma}
  \end{equation*}
  such that
  \begin{enumerate}
  \item $h(-,1) = \id$,
  \item $h(-,0) = j_{\sigma}r_{\sigma}$,
  \item $h(-,t)$ restricts to the identity on $T_{\sigma}$ for every
    $t$.
  \end{enumerate}
  So, the morphism $h$ gives an algebraic homotopy equivalence between
  $X_{\sigma}$ and $T_{\sigma}$.
\end{prp}
\begin{prf}
  We have $X_{\sigma} = \Spec k[\sigma^{\vee}\cap M]$, $T_{\sigma} =
  \Spec k[\sigma^{\bot} \cap M]$. The inclusion $T_{\sigma} \to
  X_{\sigma}$ is given by the quotient
  \begin{equation*}
    k[\sigma^{\vee}\cap M] \to k[\sigma^{\bot} \cap M],
  \end{equation*}
  which is the identity on $\sigma^{\bot}$ and sends any element $v
  \in \sigma^{\vee}$ not in $\sigma^{\bot}$ to $0 \in k[\sigma^{\bot}
  \cap M]$.

  Pick $u_0 \in \sigma$ such that $\sigma^{\bot} = \sigma^{\vee} \cap
  u_0^{\bot}$. Then define
  \begin{equation*}
    h^*\cl k[\sigma^{\vee} \cap M] \to k[\sigma^{\vee}\cap M]\tsr k[t],
  \end{equation*}
  by $h^*(v) = v \tsr t^{\angbk{u_0,v}}$ for every $v \in
  \sigma^{\vee}$.  This gives a morphism of schemes $h\cl
  X_{\sigma}\tms \bbA^1_k \to X_{\sigma}$ with the desired properties.
\end{prf}

\subsection{Morphic cohomology of an algebraic torus}

Now we compute the morphic cohomology ring of an algebraic torus. As
we will need this computation for subtori of a toric variety, it will
be useful to have a canonical description of this ring in terms of the
lattice defining the toric variety.

Let $N$ be a lattice of rank $n$, and $L_{\bbR} \subset
M_{\bbR}=N_{\bbR}^{\vee}$ a subspace of dimension $r$ generated by
vectors in the lattice $M$. Consider the rank $r$ sublattice $L =
L_{\bbR} \cap M$, and its associated torus
\begin{equation*}
  T_L = \Spec \bbC[L].
\end{equation*}

As we have seen in \ref{prp:comp}, $LH^1(\bbG_m) \simeq \bbZ[s] e$
with $e$ a generator in bidegree $(1,1)$ which corresponds, by
duality, to a radial Borel-Moore chain joining $0$ and $\oo$ in
$\bbG_m$. Now, any $v \in L$ defines a character $\chi_v\cl T_L \to
\Spec \bbC[v,v^{-1}] = \bbG_m$, so we can now define a morphism of
rings
\begin{equation}
  \varphi\cl (\wedge^* L) \tsr LHpt \to LH^*(T_L)
  \label{eq:wedgetori}
\end{equation}
by
\begin{equation*}
  \varphi(v\tsr 1) = \chi_v^*(e),
\end{equation*}
for $v \in L$, and extended in the obvious way to the exterior algebra
because $LH^*(T_L)$ is a graded commutative algebra (theorem
\ref{thm:cup}). On the left hand side of \eqref{eq:wedgetori} the
vectors in the lattice $L$ are assumed to have bigrading $(1,1)$.

\begin{thm}
  \label{thm:mortori}
  The morphism $\varphi$ is an isomorphism.
\end{thm}
\begin{prf}
  We argue by induction on the rank of $L$. The isomorphism is clear
  when $\rank L = 1$ by the computation in \ref{prp:comp}. Let $L =
  L_0 \oplus \bbZ v$. This gives a product decomposition $T_L =
  T_{L_0} \tms \bbG_m$. Now, because the Kunneth isomorphism in
  \ref{prp:kunneth} preserves the cup product, we get a commutative
  diagram
\begin{equation*}
  \xymatrix{\wedge^nL \tsr \LHpt \ar[r] \ar[d]^{\varphi} &
    \bigoplus_{p+q=n}(\wedge^p L_0 \tsr \LHpt) \tsr_{\LHpt} (\wedge^q \bbZ v \tsr
    \LHpt) \ar[d] \\
    LH^n(T_L) \ar[r] & \bigoplus_{p+q=n}LH^p(T_{L_0}) \tsr_{\LHpt} LH^q(\bbG_m)}
\end{equation*}
The upper row is an isomorphism by multilinear algebra results, while
the lower row is an isomorphism by the Kunneth isomorphism
\ref{prp:kunneth}. The right vertical map is a sum of tensor products
of $\varphi$'s corresponding to lower dimensional tori, so are
isomorphisms by induction hypothesis. We conclude then that the left
vertical map is an isomorphism.
\end{prf}

\section{Spectral sequence associated to a toric variety}
\label{sec:sseq}

Let $X(\Delta)$ be a complete toric variety of dimension $n$, $R$ a
ring (possibly graded) and $\caF^*$ a cochain complex of sheaves of
$R$-modules on $X$. As usual, the hypercohomology of $\caF^*$ will be
\begin{equation*}
  \bbH^n(X(\Delta), \caF^*) = H^n\Gamma(X(\Delta), \caI^*)
\end{equation*}
where $\caI^*$ is a K-injectives resolution $\caF^* \to \caI^*$.

In this section we will write down a spectral sequence converging to
the hypercohomology $\bbH^n(X(\Delta), \caF^*)$ whose $E_2$ page is
computable in terms of the combinatorics of the toric variety, and the
hypercohomology of $\caF^*$ on algebraic tori. The spectral sequence
comes from the identification
\begin{equation*}
  \bbH^n(X(\Delta), \caF^*) = \Ext^n(R_X, \caF^*),
\end{equation*}
and the fact that the hyper-ext can be computed resolving either
variable. We will chose to resolve the constant sheaf $R_X$ producing
a \v{C}ech-like resolution $\check{\caC}_*(\Delta, R_X) \to R_X$ from
the combinatorics of the toric variety.

A similar idea to ours, applied to singular homology and cohomology,
was previously developed in the thesis \cite{Jor:Toric}.

\subsection{Resolution associated to a fan}

Let $X(\Delta)$ be a complete toric variety defined by a fan $\Delta$,
and let $R$ be a commutative ring.

\begin{dfn}
  \label{dfn:cechmod}
  Let $\check{C}_k(\Delta, R)$ for $k\geq 0$ be the sequence of free
  $R$-modules
  \begin{equation*}
    \check{C}_k(\Delta, R) = \bigoplus_{\sigma\in \Delta^{(k)}} R[\sigma],
  \end{equation*}
  with basis $\Delta^{(k)}$, together with $R$-morphisms $d_k\cl
  \check{C}_k(\Delta, R) \to \check{C}_{k-1}(\Delta, R)$ defined by
  \begin{equation*}
    d_k([\tau])_{\sigma} = \epsilon(\tau, \sigma),
  \end{equation*}
  where $\epsilon(\tau,\sigma) = \pm 1$ according to whether the
  orientation induced by $\sigma$ on $\tau$ coincides or not with the
  orientation in $\tau$ (remember, we assume a fixed choice of
  orientation on every cone).
\end{dfn}

\begin{dfn}
  \label{dfn:checksh}
  Let $\check{\caC}_k(\Delta, R_X)$ for $k\geq 0$ be the sequence of
  sheaves of $R$-modules on $X(\Delta)$ given by
  \begin{equation*}
    \check{\caC}_k(\Delta, R_X) = \bigoplus_{\sigma \in
      \Delta^{(k)}} i_{\sigma !}i_{\sigma}^* R_X,
  \end{equation*}
  where $R_X$ is the constant sheaf on $X(\Delta)$ and $i_{\sigma} \cl
  X_{\sigma} \to X(\Delta)$ is the inclusion of the open $X_{\sigma}$.

  Moreover, we define, a sequence of morphisms $d_k\cl
  \check{\caC}_k(\Delta, R_X) \to \check{\caC}_{k-1}(\Delta, R_X)$
  given by
  \begin{equation*}
    d_k = \bigoplus_{\substack{\sigma \in \Delta^{(k-1)}\\ \tau \in
        \Delta^{(k)}\\ \tau \leq \sigma}} \epsilon(\tau,\sigma)\iota_{\tau,\sigma}
  \end{equation*}
  where $\iota_{\tau,\sigma}\cl i_{\tau !} i_{\tau}^*R_{X} \to
  i_{\sigma!}  i_{\sigma}^* R_{X}$ is the natural inclusion of sheaves
  inducing the identity on the nonzero fibers.
\end{dfn}

\begin{dfn}
  Given a fan $\Delta$ and a cone $\sigma \in \Delta$ of codimension
  $k$, we define a fan $\Delta_{\sigma}$ defined on the lattice
  $N/(\bbR\sigma\cap N)$ of dimension $k$, whose cones are the
  projection of cones in $\Delta$ having $\sigma$ as a face.
\end{dfn}
\begin{rmk}
  Notice that the cones in $\Delta_{\sigma}$ correspond bijectively
  with the cones $\tau \in \Delta$ having $\sigma$ as a face.
\end{rmk}

Let $x \in X(\Delta)$ be a point. We denote by $\sigma(x)$ the unique
cone in $\Delta$ such that $x \in T_{\sigma(x)}$.

\begin{rmk}
\label{rmk:fibercech}
Observe that 
\begin{equation*}
  (i_{\tau !}i_{\tau}^*R_{X})_x = \begin{cases}
    R & \text{ if $\sigma(x) \leq \tau$,}\\
    0 & \text{ otherwise.}
  \end{cases}
\end{equation*} 
In other words, the fiber $(i_{\tau !}i_{\tau}^*R_{X})_x$ is nonzero
exactly for the cones $\tau \in \Delta_{\sigma(x)}$.
\end{rmk}

This allows us to define a morphism of $R$-modules
\begin{equation*}
  f_{k,x}\cl \check{\caC}_k(\Delta, R_X)_x \to \check{C}_k(\Delta_{\sigma(x)}, R)
\end{equation*}
which is the identity on every nonzero summand, as the summands on the
right correspond exactly to the nonzero summands on the left.

\begin{prp}
  \label{prp:fibers}
  Let $X(\Delta)$ be a complete toric variety associated to a fan $\Delta$.
  \begin{enumerate}
  \item The sequence of $R$-modules $\check{C}_k(\Delta, R)$ together
    with the differentials $d_k$ form a chain complex of $R$-modules,
    that is,
    \begin{equation*}
      d_{k-1}d_{k} = 0.
    \end{equation*}

  \item The diagram
    \begin{equation*}
      \xymatrix{\check{\caC}_k(\Delta, R_X)_x \ar[r]^{f_{k,x}} \ar[d]^{d_k} &
        \check{C}_k(\Delta_{\sigma(x)}, R) \ar[d]^{d_k}\\
        \check{\caC}_{k-1}(\Delta, R_X)_x \ar[r]^{f_{k-1,x}} &
        \check{C}_{k-1}(\Delta_{\sigma(x)}, R)}
    \end{equation*}
    is commutative, and the rows are isomorphisms.

  \item The sequence of sheaves of $R$-modules $\check{\caC}_k(\Delta,
    R_X)$ together with the differentials $d_k$ form a chain complex
    of sheaves of $R$-modules.
  \end{enumerate}
\end{prp}
\begin{prf}
  1) Let $[\tau]$ be an element of the basis of $\check{C}_k(\Delta,
  R)$. Then we have
\begin{eqnarray*}
  d_{k-1}d_k([\tau])_{\eta} &=& \sum_{\eta < \sigma < \tau} \epsilon(\eta, \sigma)
  \epsilon(\sigma,\tau) 
\end{eqnarray*}
and either there is no $\sigma$ in between $\eta$ and $\tau$ or there
are exactly two of them, giving opposite signs.

2) Follows easily from the definition and remark \ref{rmk:fibercech}.

3) From 1 and 2 together we see that $d_{k-1,x}d_{k,x}\cl
\check{\caC}_k(\Delta, R_X)_x \to \check{\caC}_{k-2}(\Delta, R_X)_x$
is zero. Then as $(d_{k-1}d_k)_x = d_{k-1,x}d_{k,x}$ and a morphism of
sheaves which is zero on the fibers is the zero morphism, we are done.
\end{prf}

\begin{dfn}
  \label{dfn:aug}
  Let $a\cl \check{C}_*(\Delta, R) \to R$ be the augmentation morphism
  given by
  \begin{equation}
    a(x) = \sum_{\tau \in \Delta^{(0)}} x_{\tau},
  \end{equation}
  where $x = \sum x_{\tau}[\tau] \in \check{C}_0(\Delta, R)$.

  In a similar way, let $a\cl \check{\caC}_*(\Delta, R_X) \to R_X$ be
  the augmentation morphism induced by the morphisms $i_{\sigma
    !}i_{\sigma}^* R_{X} \to R_{X(\Delta)}$.
\end{dfn}

We will prove that $a\cl \check{\caC}_*(\Delta, R_X) \to R_X$ is a
quasi-isomorphism. To do so, we will relate the fiber complexes
$\check{\caC}_*(\Delta, R_X)_x$ with the cellular homology complex of
a cellular decomposition on the $n$-dimensional ball.

Let
\begin{eqnarray}
  B &=& \set{x \in N_{\bbR}\mid \norm{x} \leq 1},\\
  S &=& \set{x \in N_{\bbR} \mid \norm{x} = 1}.
\end{eqnarray}
\begin{dfn}
  \label{dfn:cellsphere}
  For every non-zero cone $\sigma \in \Delta$ we define the subset
  $e_{\sigma}\subset S$ as follows
  \begin{equation*}
    e_{\sigma} = \sigma \cap S.
  \end{equation*}
\end{dfn}

\begin{prp}
  \label{prp:cellsphere}
  The subsets $e_{\sigma} \subset S$ are homeomorphic to closed balls
  of dimension $\dim \sigma - 1$. Together form a cellular
  decomposition of the sphere $S$. The sphere together with this
  decomposition will be denoted by $S_{\Delta}$.
\end{prp}
\begin{proof}
  As $e_{\sigma}$ are the intersection of a strictly convex cone with
  the unit sphere $S_{\Delta}$ it is clear that are homeomorphic to
  balls of dimension $\dim \sigma - 1$. They cover all $S$ because the
  fan is complete. Finally, the boundary of the cell $e_{\sigma}$ is
  formed by the cells associated to the faces of $\sigma$, and so
  belong to the lower dimensional skeleton.
\end{proof}

\begin{dfn}
  \label{dfn:simpcpx}
  To any complete fan $\Delta$ we associate an abstract simplicial
  complex $K(\Delta)$ as follows:
  \begin{enumerate}
    \item The vertices in $K(\Delta)$ correspond to the cones in
      $\Delta$.
    \item The $k$-simplexes in $K(\Delta)$ are the sets of vertices
      belonging to flags in $\Delta$ of length $k$, that is, sequences
      of strictly included cones
      \begin{equation*}
        \tau_0 < \tau_1 < \cdots < \tau_k.
      \end{equation*}
  \end{enumerate}  
\end{dfn}
\begin{rmk}
  The simplicial complex $K(\Delta)$ is related to the barycentric
  subdivision of the fan $\Delta$. However, they are not exactly the
  same thing because of the cone $0 \in \Delta$. If we had not admited
  the cone $0$ in the definition of $K(\Delta)$ we would have obtained
  a combinatorial model of the barycentric subdivision of the fan
  $\Delta$.
\end{rmk}

For every $1$-dimensional cone in $\tau \in \Delta^{(n-1)}$ let
$u_{\tau} \in N_{\bbR}$ be the unique unit vector generating it. Then,
for any non-zero cone $\sigma \in \Delta$, let $v_{\sigma}$ be the
vector
\begin{equation*}
  v_{\sigma} = \sum_{\substack{\tau \in \Delta^{(n-1)}\\ \tau \leq
      \sigma }} u_{\tau}.
\end{equation*}
\begin{dfn}
  \label{dfn:cellball}
  For every $k$-simplex $(\tau_0,\ldots,\tau_k) \in K(\Delta)$ given
  by a flag of cones $\tau_0 < \cdots < \tau_k$, we define a subset
  $d_{(\tau_0,\ldots,\tau_k)} \subset B$ as follows,
  \begin{equation}
    d_{(\tau_0,\ldots,\tau_k)} =  \begin{cases}
      \set{0} & \text{ if $\tau_0 = 0$ and $k = 0$,}\\
      \bbR_{\geq
        0} \angbk{v_{\tau_1}, \ldots, v_{\tau_k}}\cap B & \text{ if
        $\tau_0 = 0$ and $k > 0$,}\\
      \bbR_{\geq
        0} \angbk{v_{\tau_0}, \ldots, v_{\tau_k}}\cap S & \text{ if $\tau_0 \neq 0$.}
    \end{cases}
  \end{equation}
\end{dfn}
\begin{prp}
  \label{prp:cellball}
  The subsets $d_{\tau_0,\ldots,\tau_k} \subset B$ are homeomorphic to
  closed balls of dimension $k$. Together form a cellular
  decomposition of the ball $B$. The ball together with this
  decomposition will be denoted by $B_{K(\Delta)}$.
\end{prp}
\begin{prf}
  Let $(\tau_0,\ldots,\tau_k) \in K(\Delta)$. Because the vectors
  $v_{\tau_i}$ all belong to the cone $\tau_k$, the subsets
  $\bbR_{\geq 0} \angbk{v_{\tau_0}, \ldots, v_{\tau_k}}$ are strongly
  convex cones. Then arguing as in \ref{prp:cellsphere} for every case
  in definition \ref{dfn:cellball} we prove that
  $d_{\tau_0,\ldots,\tau_k}$ are balls. As for the statment regarding
  its dimension, it follows from the linear independence of the
  $v_{\tau_i}$ for any flag $0 \neq \tau_0 < \cdots < \tau_k$.

  Finally, observe that the boundary of a cell
  $d_{\tau_0,\ldots,\tau_k}$ is formed by the cells resulting from
  removing one cone in the flag, all of lower dimension. This proves
  that the cells $d_{\tau_0,\ldots,\tau_k}$ give a cellular
  decomposition of the ball $B$.
\end{prf}

\begin{dfn}
  \label{dfn:celldual}
  For every cone $\sigma \in \Delta$, let $e_{\sigma}^{\vee} \subset
  B$ be the subset defined by
  \begin{equation}
    e_{\sigma}^{\vee} = \bigcup_{\substack{k \geq 0\\(\tau_0,\ldots,\tau_k) \in
        K(\Delta)\\\sigma \leq \tau_0}} d_{\tau_0,\ldots,\tau_k}
  \end{equation}
\end{dfn}
\begin{prp}
  \label{prp:celldual}
  The subsets $e_{\sigma}^{\vee} \subset B$ are homeomorphic to closed
  balls of dimension $\codim \sigma$. Together form a cellular
  decomposition of the ball $B$. The ball together with this
  decomposition will be denoted by $B_{\Delta}^{\vee}$.
\end{prp}
\begin{prf}
  Observe that $e_{\sigma}^{\vee}$ is a geometric realization of a
  subcomplex of $K(\Delta)$. Moreover, this subcomplex happens to be
  isomorphic to $K(\Delta_{\sigma})$ (follows directly from the
  definitions). Now, applying proposition \ref{prp:cellball} to the
  simplicial complex $K(\Delta_{\sigma})$ we get a cellular
  decomposition of a $(\codim \sigma)$-dimensional ball
  $B_{\Delta_{\sigma}}$ which realizes the simplicial complex
  $K(\Delta_{\sigma})$. We have two geometric realizations of
  isomorphic simplicial complexes, so they must be isomorphic, and we
  conclude that $e_{\sigma}^{\vee}$ is homeomorphic to a $(\codim
  \sigma)$-dimensional ball.

  It is clear that the cells $e_{\sigma}^{\vee}$ cover all the ball
  $B$, and that they are attached properly, that is the boundary of
  every cell is contained in a lower dimensional skeleton.
\end{prf}
\begin{rmk}
  The cellular decomposition $B_{\Delta}^{\vee}$ is dual to the
  cellular decomposition of the sphere $S_{\Delta}$ together with an
  extra $n$-cell for the interior of the ball.
\end{rmk}

\begin{prp}
  \label{prp:cpxiso}
  There is a canonical isomorphism of chain complexes
  \begin{equation*}
    \check{C}_*(\Delta, R) \simeq C^{\mrm{cell}}_*(B_{\Delta}^{\vee}, R).
  \end{equation*}
\end{prp}
\begin{prf}
  There is a canonical isomorphism between the $R$-modules
  $\check{C}_k(\Delta, R)$ and $C^{\mrm{cell}}_k(B_{\Delta}^{\vee},
  R)$, as both are generated by the cones in codimension $k$. Note
  that the unique $n$-cell in $B_{\Delta}$ corresponds to the zero
  cone in $\Delta$. 

  It only remains to check that the differentials in
  $\check{C}_*(\Delta, R)$ coincide with the cellular ones. Note that
  the attaching maps $f_{\tau}\cl \partial e_{\tau}^{\vee} \to
  \mrm{Sk}_{\codim \tau - 1}B_{\Delta}^{\vee}$ are homeomorphisms with
  the image. So, for any lower dimensional cell $e_{\sigma}^{\vee}$ on
  the boundary of $e_{\tau}^{\vee}$, the corresponding matrix element
  in the cellular differential is a sign, according to the relative
  orientation of the cells $e_{\sigma}^{\vee}$ and
  $e_{\tau}^{\vee}$. This is exactly the differential in
  $\check{C}_*(\Delta, R)$.
\end{prf}

\begin{cor}
  \label{cor:toricres}
  The augmented complex of sheaves
  \begin{equation*}
    \check{\caC}_*(\Delta, R_X) \to R_X
  \end{equation*}
  is exact, that is, $\check{\caC}_*(\Delta, R_X)$ is a resolution of
  the constant sheaf $R_X$.
\end{cor}
\begin{prf}
  From proposition \ref{prp:fibers} we know that the fiber of the
  complex $\check{\caC}_*(\Delta, R_X) \to R_X$ at a point $x \in
  X(\Delta)$ is identified with $\check{C}(\Delta, R) \to R$. Now by
  proposition \ref{prp:cpxiso} we know that the complex
  $\check{C}(\Delta, R)$ is isomorphic to the cellular complex of the
  cellular decomposition of the ball $B_{\Delta}^{\vee}$, so its
  homology is
  \begin{equation*}
    H_k\check{C}(\Delta, R) \simeq H_k C^{\mrm{cell}}(B_{\Delta}, R)
    = \begin{cases} R & \text{for $k = 0$,}\\
      0 & \text{for $k > 0$.} \end{cases}
  \end{equation*}
  Now, as the augmentation $a\cl \check{C}_*(\Delta, R) \to R$ is
  given by the sum of the elements on degree zero, we conclude that it
  is a quasi-isomorphism.
\end{prf}

\subsection{The spectral sequence}

Let $X(\Delta)$ be a toric variety and $\caF^*$ be a complex of
sheaves on $X$. We describe a spectral sequence converging to the
hypercohomology $\bbH^n(X(\Delta), \caF^*)$.

\begin{dfn}
  A complex of sheaves $\caF^*$ is said to have \emph{homotopy
    invariant cohomology} if for every variety $X$ the projection
  $p\cl X\tms \bbA^1 \to X$ induces isomorphisms in hypercohomology
  \begin{equation*}
    \bbH^n(X, \caF^*) \weq \bbH^n(X\tms \bbA^1, \caF^*).
  \end{equation*}
\end{dfn}

\begin{rmk}
  The complex of sheaves $\caM^*$ defining morphic cohomology has
  homotopy invariant cohomology by \ref{thm:morhiHI}.
\end{rmk}

\begin{thm}
  \label{thm:sseq}
  Let $X(\Delta)$ be a toric variety associated to a fan $\Delta$ and
  $\caF^*$ a bounded above cochain complex of sheaves. There is a
  convergent spectral sequence
  \begin{equation}
    E_1^{r,s} =\Ext^r(\check{\caC}_s(\Delta, \bbZ_X), \caF^*)
    \Longrightarrow \bbH^{r+s}(X(\Delta), \caF^*).
    \label{eq:sseqTor}
  \end{equation}

  Moreover, if $\caF^*$ has homotopy invariant cohomology,
  \begin{equation}
    E_1^{r,s} = \bigoplus_{\sigma \in \Delta^{(s)}} \bbH^r(T_{\sigma},
    \caF^*),
    \label{eq:sseqHI}
  \end{equation}
  and the differentials on the first page $d_1\cl E_1^{r,s} \to
  E_1^{r,s+1}$ are given by
  \begin{equation}
    d_1 = \sum_{\substack{\sigma \in \Delta^{(s)}\\\tau \in
        \Delta^{(s+1)} \\ \tau \leq \sigma}} \epsilon(\tau,\sigma)
      r_{\tau,\sigma}^*
    \label{eq:sseqHIdiff}
  \end{equation}
  where
  \begin{equation*}
    r_{\tau,\sigma}\cl T_{\tau} = \Spec \bbC[\tau^{\bot}] \to T_{\sigma} = \Spec \bbC[\sigma^{\bot}],
  \end{equation*}
  are the morphisms induced by the natural inclusion $\sigma^{\bot}
  \to \tau^{\bot}$.
\end{thm}
\begin{prf}
  Let $\caF^* \to \caI^{*}$ be a K-injective resolution of $\caF^*$
  (see remark \ref{rmk:resolutions}). Let $\check{\caC}_*(\Delta,
  \bbZ_X) \to \bbZ_X$ be the resolution of the constant sheaf $\bbZ_X$
  from corollary \ref{cor:toricres}. We build a double complex
  \begin{equation*}
    C^{r,s} = \Hom(\check{\caC}_s(\Delta, \bbZ_X), \caI^r),
  \end{equation*}
  with the induced differentials (going in the increasing direction of
  $r$ and $s$). The homology of this double complex in the $r$
  direction is $\Ext^r(\check{\caC}_s(\Delta, \bbZ_X), \caF^*)$, giving the
  spectral sequence
  \begin{equation*}
    E_1^{r,s} =\Ext^r(\check{\caC}_s(\Delta, \bbZ_X), \caF^*) \Longrightarrow
    \bbH^{r+s}(X(\Delta), \caF^*).
  \end{equation*}

  As for the convergence. The complex of sheaves $\caI^*$ is bounded
  above, and the schemes $X_{\sigma}$ have finite cohomological
  dimension. Using the hypercohomology spectral sequence we conclude
  that $\bbH^k(X_{\sigma}, \caF^*)$ is bounded above. In other words,
  the first page is bounded above in the $r$ direction. It is clearly
  bounded (from both sides) in the $s$ direction, and this is enough
  to establish the convergence.

  If $\caF^*$ is homotopy invariant, as the immersion $T_{\sigma} \to
  X_{\sigma}$ are algebraic homotopy equivalences we get the
  expression \eqref{eq:sseqHI}. 

  Finally, the differentials on the first page are induced by the
  $s$-differentials in the double complex $C^{r,s}$, which are given
  by the formula
  \begin{equation*}
    d_1 = \sum_{\substack{\sigma \in \Delta^{(s)}\\\tau \in
        \Delta^{(s+1)} \\ \tau \leq \sigma}} \epsilon(\tau,\sigma) i_{\tau,\sigma}^*
  \end{equation*}
  where $i_{\tau,\sigma}\cl X_{\tau} \to X_{\sigma}$ is the
  inclusion. The formula \eqref{eq:sseqHIdiff} follows from the
  equation
  \begin{equation*}
    r_{\tau,\sigma}^* = j_{\tau}^* i_{\tau,\sigma}^* r_{\sigma}^*,  
  \end{equation*}
  and the fact that $j_{\tau}^*$ and $r_{\sigma}^*$ are mutually
  inverse isomorphisms giving the identification $\bbH^r(T_{\sigma},
  \caF^*) \simeq \bbH^r(X_{\sigma}, \caF^*)$.
\end{prf}

\subsection{Some applications}
We have a rather explicit description of the first page and
differentials of the spectral sequence in \ref{thm:sseq}. Together
with the computation \ref{thm:mortori} of the morphic cohomology of a
torus we can make it still more explicit.
\begin{cor}
  \label{cor:sseqmor}
  Let $\caF^* = \caM^*$. Then the first page of the spectral sequence
  in \ref{thm:sseq} is
  \begin{equation}
    E_1^{r,s} \simeq \bigoplus_{\sigma \in \Delta^{(s)}} \wedge^r(\sigma^{\bot}
    \cap M) \tsr \LHpt
    \label{eq:sseqmor}
  \end{equation}
  and the differentials $d^s_1\cl E_1^{r,s} \to E_1^{r,s+1}$ are given
  by
  \begin{equation*}
    d^s_1(\sum_{\sigma \in \Delta^{(s)}} x_{\sigma}
    v_{1,\sigma}\wedge\cdots \wedge
    v_{r,\sigma}) = \sum_{\sigma\in \Delta^{(s)}} x_{\sigma} \sum_{\substack{\tau \in \Delta^{(s+1)}\\          \tau \leq \sigma}}
    \epsilon(\tau, \sigma) v_{1,\sigma}\wedge \cdots
    \wedge v_{r,\sigma}
  \end{equation*}
\end{cor}
\begin{prf}
  This is a straight forward combination of theorem \ref{thm:sseq} and
  the computation \ref{thm:mortori}.
\end{prf}

Finally, using an idea from \cite{Jor:Toric} which can be traced back
to \cite{Tot:ChowLinear} that this spectral sequence degenerates
rationally.
\begin{thm}
  \label{thm:sseqrat}
  The spectral sequence in \ref{cor:sseqmor} degenerates when tensored
  with $\bbQ$.
\end{thm}
\begin{prf}
  Let $\LHpt_{\bbQ} = \LHpt\tsr_\bbZ \bbQ$. The toric variety
  $X(\Delta)$ admits an $\bbN$-action. Let $m \in \bbN$, then $[m]\cl
  X(\Delta) \to X(\Delta)$ is the morphism which on the open sets
  $X_{\sigma} = \Spec \bbC[\sigma^{\vee}\cap M]$ is defined through
  the ring homomorphism $[m]^*\cl \bbC[\sigma^{\vee}\cap M] \to
  \bbC[\sigma^{\vee}\cap M]$ given by $v \mapsto mv$ (see
  \cite{Tot:ChowLinear} for details).

  The $\bbN$-action on $X(\Delta)$ induces an $\bbN$-action on the
  spectral sequence \ref{eq:sseqmor}. As the rational morphic
  cohomology of a torus $T_L = \Spec \bbC[L]$ is
  \begin{equation*}
    LH^n(T_L)_{\bbQ} = \wedge^n L \tsr \LHpt_{\bbQ},
  \end{equation*}
  The $\bbN$-action on on the page $E_1^{r,s}$ is just multiplication
  by $m^r$. As the next pages $E_k^{r,s}$ of the spectral sequence are
  subquotients of $E_1^{*,*}$, the action on those pages is also given
  by $m^r$. On the other hand the differentials go $d_k\cl E_k^{r,s}
  \to E_k^{r+1-k,s+k}$ and the $\bbN$-action commutes with them, so
  \begin{equation*}
    m^rd_k(x) = d_k(m^rx) = d_k([m]x) = [m]d_k(x) = m^{r+1-k}d_k(x),
  \end{equation*}
  where $x\in E_k^{r,s}$. Rationally, this forces $d_k(x) = 0$ when $k
  \geq 2$.
\end{prf}


Now let $\varepsilon \cl \ctTop \to \ctqProj_\bbC$ be the morphism of
sites, with the usual topology in $\ctTop$ and the Zariski topology on
$\ctqProj_{\bbC}$. Let $\mbf{R}\varepsilon_* \bbZ$ be the derived
push-forward of the constant sheaf $\bbZ$ on $\ctTop$ to the Zariski
site $\ctqProj_{\bbC}$. There is a natural map
\begin{equation}
  \caM^*(q) \to \mbf{R}\varepsilon_* \bbZ
  \label{eq:suslin}
\end{equation}
which, on smooth varieties, factors as
\begin{equation}
  \caM^*(q) \to \tau_{\leq q}\mbf{R}\varepsilon_* \bbZ.
  \label{eq:suslintr}
\end{equation}
See \cite{FriHaWa:TechKth} for details.

There is the following conjecture, a morphic analogue of the
Beilinson-Lichtenbaum conjecture in the motivic world.
\begin{con}[Suslin]
  The comparison morphism \eqref{eq:suslin} above is a
  quasi-isomorphism on smooth varieties.
\end{con}

This conjecture is proved for the class of smooth linear varieties
(which include smooth toric varieties) in \cite{FriHaWa:TechKth}
theorem 7.14.

The spectral sequence \ref{thm:sseq} has the following corollary.
\begin{cor}
  The Suslin conjecture holds for projective toric varieties (not
  necessarily smooth).
\end{cor}
\begin{prf}
  First of all, we have to check that $\caM^*(q)\mid_{X(\Delta)}$ is
  exact above degree $q$, in order to have the factorization
  \eqref{eq:suslintr}.
  This is a local statment on $X(\Delta)$, so we can restrict to an
  open $X_{\sigma}$. Now the inclusion $T_{\sigma} \to X_{\sigma}$ is
  an algebraic homotopy equivalence, and they induce isomorphisms on
  hypercohomology
  \begin{equation*}
    \bbH^n(X_{\sigma}, \caM^*(q)|_X) \weq \bbH^n(T_{\sigma},
    \caM^*(q)|_{T_{\sigma}}),
  \end{equation*}
  so the natural map $\caM^*(q)|_X \to \mbf{R} j_{\sigma *}
  j_{\sigma}^* \caM^*(q)|_X$ is a quasi-isomorphism. As $T_{\sigma}$
  is smooth, its cohomology vanishes above $q$, and we have the
  desired factorization.
  
  Now, $\tau_{\leq q}\mbf{R}\varepsilon_*\bbZ$ is a has homotopy
  invariant cohomology. Then we can apply theorem \ref{thm:sseq} to it
  and get a spectral sequence converging to $\bbH^n(X(\Delta),
  \tau_{\leq q}\mbf{R}\varepsilon_*\bbZ)$. Moreover, the comparison
  map \eqref{eq:suslintr} gives a morphism of spectral sequences
  \begin{equation*}
    \xymatrix{E_1^{r,s} = \bigoplus_{\sigma \in \Delta^{(s)}}
      \bbH^r(T_{\sigma}, \caM^*(q)) \ar@2{->}[r] \ar[d]&
      \bbH^{r+s}(X(\Delta), \caM^*(q)) \ar[d]\\
      E_1^{'r,s} = \bigoplus_{\sigma \in \Delta^{(s)}}
      \bbH^r(T_{\sigma}, \tau_{\leq q}\mbf{R}\varepsilon_*\bbZ) \ar@2{->}[r]
      &\bbH^{r+s}(X(\Delta), \tau_{\leq q}\mbf{R}\varepsilon_*\bbZ)}
  \end{equation*}
  which is an isomorphism on the first page by the computation
  \ref{thm:mortori}, so it gives an isomorphism on the abutement, as
  desired.
\end{prf}



\bibliographystyle{alphaurl}

\bibliography{toric}

\makeauthors

\end{document}